\documentclass[12pt,oneside,reqno]{amsart}
\usepackage{graphicx}
\usepackage{mathrsfs}
\usepackage{stmaryrd}
\usepackage{amsfonts}
\usepackage{enumerate,amsmath,amssymb,amsthm}
\newcommand{\be}{\begin{eqnarray}}
\newcommand{\ee}{\end{eqnarray}}
\newcommand{\ce}{\begin{eqnarray*}}
\newcommand{\de}{\end{eqnarray*}}
\newtheorem{theorem}{Theorem}[section]
\newtheorem{lemma}[theorem]{Lemma}
\newtheorem{remark}[theorem]{Remark}
\newtheorem{definition}[theorem]{Definition}
\newtheorem{proposition}[theorem]{Proposition}

\newtheorem{corollary}[theorem]{Corollary}

\def\a{\alpha}
\def\om{\omega}

\def\b{\beta}
\def\p{\partial}

\def\[{{\Big[}}
\def\]{{\Big]}}
\def\<{{\langle}}
\def\>{{\rangle}}
\def\({{\Big(}}
\def\){{\Big)}}

\def\dif{{\mathord{{\rm d}}}}

\def\div{{\mathord{{\rm div}}}}

\def\u{\mathord{{\bf u}}}
\def\f{\mathord{{\bf f}}}
\def\v{\mathord{{\bf v}}}
\def\e{\mathord{{\bf e}}}
\def\y{\mathord{{\bf y}}}
\def\b{\mathord{{\bf b}}}
\def\w{\mathord{{\bf w}}}

\def\no{\nonumber}
\def\bt{\begin{theorem}}
\def\et{\end{theorem}}
\def\bl{\begin{lemma}}
\def\el{\end{lemma}}
\def\br{\begin{remark}}
\def\er{\end{remark}}
\def\bd{\begin{definition}}
\def\ed{\end{definition}}
\def\bp{\begin{proposition}}
\def\ep{\end{proposition}}
\def\bc{\begin{corollary}}
\def\ec{\end{corollary}}

\def\cO{{\mathcal O}}

\def\cV{{\mathcal V}}

\def\mH{{\mathbb H}}

\def\mN{{\mathbb N}}

\def\mR{{\mathbb R}}

\def\geq{\geqslant}
\def\leq{\leqslant}

\def\sP{{\mathscr P}}

\def\bR{{\mathbf R}}

\def\lb{{\llbracket}}
\def\rb{{\rrbracket}}

\allowdisplaybreaks

\begin{document}

\title{Tamed 3D Navier-Stokes Equation: Existence, Uniqueness and Regularity}

\date{}
\author{Michael R\"ockner, Xicheng Zhang }

\dedicatory{
Fakult\"at f\"ur Mathematik,
Universit\"at Bielefeld\\
Postfach 100131,
D-33501 Bielefeld, Germany\\
 M. R\"ockner: roeckner@math.uni-bielefeld.de\\
X. Zhang: xzhang@math.uni-bielefeld.de
 }

\begin{abstract}

In this paper, we prove the existence and uniqueness of a smooth solution to a
tamed 3D Navier-Stokes equation in the whole space. In particular, if there exists a bounded
smooth solution to the classical 3D Navier-Stokes equation, then this solution satisfies our tamed
equation. Moreover, using this renomalized equation we can give a new construction for a
suitable weak solution of the classical 3D Navier-Stokes equation introduced in \cite{Sc} and \cite{Ca-Ko-Ni}.
\end{abstract}

\maketitle \rm
\section{Introduction and Main Results}

Let $\u(t,x)=(u^1(t,x),u^2(t,x),u^3(t,x))$ be a row vector valued function
on $[0,\infty)\times\mR^3$. The following notations will be used throughout this paper:
$$
|\u|^2:=\sum^3_{i=1}|u^i|^2,\quad
\p_tu^j:=\frac{\p u^j}{\p t},\quad \p_iu^j:=\frac{\p u^j}{\p x_i}, \quad i,j=1,2,3
$$
$$
\nabla u^j(t):=(\p_1u^j(t),\p_2u^j(t),\p_3u^j(t)),
\quad \Delta u^j(t):=\sum_{i=1}^3\p^2_iu^j(t),\quad
j=1,2,3
$$
$$
\div(\u(t)):=\sum_{i=1}^3 \p_i u^i(t),\quad (\u(t)\cdot\nabla)\u(t):=\sum_{i=1}^3u^i(t)\p_i \u(t).
$$

Consider the following Navier-Stokes equation in $\mR^3$
\be
\p_t\u(t)=\nu\Delta \u(t)-(\u(t)\cdot\nabla)\u(t)+\nabla p(t)+\f(t),\label{NS1}
\ee
subject to the incompressibility condition:
\be
\div (\u(t))=0,\label{NS3}
\ee
and initial conditions:
\be
\u(0,x)=\u_0(x),\quad \div(\u_0)=0.\label{NS2}
\ee
Here $\u(t,x)$ represents the velocity field, $\nu>0$ is the viscosity constant,
the pressure $p(t,x)$ is an unknown scalar function,
and the external force $\f$ is a known vector valued function.

The concepts of weak solutions of (\ref{NS1})-(\ref{NS2}) and their regularities were already introduced
in the fundamental paper of Leray \cite{Le}. Pioneering work  of Leray \cite{Le} and Hopf \cite{Ho}
showed that for any $\u_0\in L^2(\mR^3)$ and $\f\in L^2([0,T]\times\mR^3;\mR^3)$,
there exist functions $\u$ and $p$ such that

\begin{enumerate}[(i)]
\item $\u\in L^\infty([0,T]; L^2(\mR^3;\mR^3))\cap L^2([0,T]; \mH^1)$;
\item the equation (\ref{NS1}) holds for $\u, p$ in the sense of distribution;
\item for almost all $t\in[0,T]$, the following energy inequality holds
\ce
\|\u(t)\|^2_{L^2}+2\nu\int^t_0\|\u(s)\|^2_{\mH^1}\dif s\leq \|\u_0\|^2_{L^2}
+2\int^t_0\<\f(s),\u(s)\>_{L^2}\dif s;
\de
\item $\lim_{t\downarrow 0}\|\u(t)-\u_0\|_{L^2}=0$.
\end{enumerate}
Here $\mH^1$ stands for the Sobolev space of order one in $L^2$ of divergence free vector fields
on $\mR^3$.

It is well known that such $\u$ is unique for the $2$D Navier-Stokes equation.
However, in the case of three dimensions, the uniqueness is only proved for small
initial data $\u_0$ and $\f$
or for large enough $\nu$. On the other hand, for any $\nu$ and $\u_0,\f$, Fabes-Jones-Rivi\`ere
\cite{Fa-Jo-Ri}(see also \cite{He}) proved that
for some small $T_*$(depending on $\nu$ and $\u_0,\f$), there exists a unique smooth solution to
(\ref{NS1})-(\ref{NS2}) on the time interval $[0,T_*)$.
To compensate for the  non-linear term $(\u\cdot\nabla)\u$, let us consider the following tamed
Navier-Stokes equation in $\mR^3$:
\be
\p_t\u(t)=\nu\Delta \u(t)-(\u(t)\cdot\nabla)\u(t)+\nabla p(t)-g_N(|\u(t)|^2)\u(t)+\f(t),\label{NS}
\ee
subject to the incompressibility condition:
\be
\div (\u(t))=0,
\ee
and initial conditions:
\be
\u(0,x)=\u_0(x),\quad \div(\u_0)=0,\label{NS0}
\ee
where $N>0$ and $g_N:\mR^+\mapsto\mR^+$ is a smooth function such that
\be
\label{Con}\left\{
\begin{array}{ll}
g_N(r):=0, \quad r\in[0,N], \\
g_N(r):=\frac{r-N-\frac{1}{2}}{\nu},\quad  r\geq N+1,\\
0\leq g'_N(r)\leq C_\nu,\quad  r\geq 0,\\
|g^{(k)}_N(r)|\leq C_{\nu,k},\quad  r\geq 0, k\in\mN,
\end{array}
\right.
\ee
that is, $g_N$ qualitatively looks as follows for $\nu<1$:
\begin{center}
\includegraphics[width=2.9in]{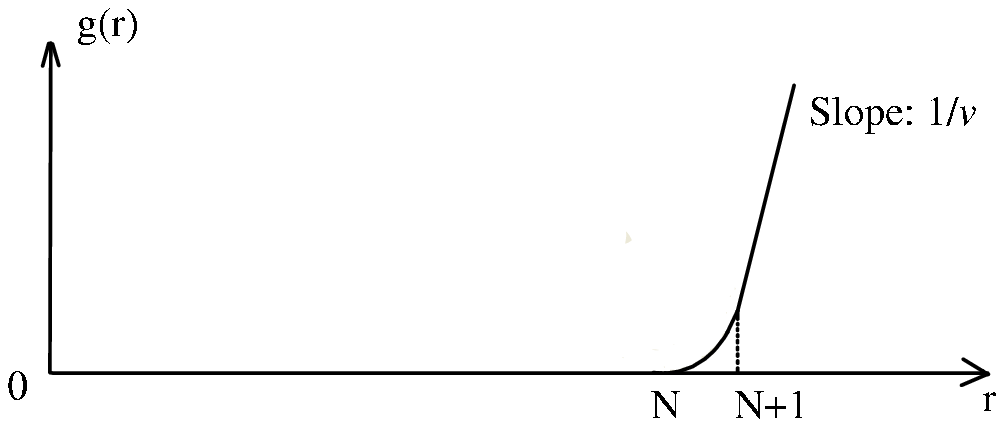}
\end{center}
\begin{center}
{\tiny  Figure: taming function: $r\mapsto g_N(r)$}
\end{center}

Here and below we shall use the following convention:  The letter $C$ with subscripts will
denote a constant depending only on its subscripts. The letter $C$ without subscripts
will denote an absolute constant, i.e., its value
does not depend on any data. All the constants may have
different values in different places.

Let us say some words about the taming function $g_N$. From the construction of
$g_N$, one  sees that if $\u(t,x)$ is a bounded smooth solution of (\ref{NS1})-(\ref{NS2})
(say bounded by $\sqrt{N}$),
then $\u(t,x)$ also satisfies (\ref{NS})-(\ref{NS0}). Intuitively,
when the velocity of the fluid is larger than $\sqrt{N}$, the dissipative term
$g_N(|\u(t)|^2)\u(t)$(regarded as some extra force) will enter into the equation and
restrain the flux of the liquid. In this sense, the value of $N$ plays the role of a valve.
As stated in \cite{Con}, the velocities observed in turbulent flows on earth are bounded.
If one accepts this as a
physical assumption, and if the system (\ref{NS1})-(\ref{NS2})
precisely describes the physical phenomenon, the solution of (\ref{NS1})-(\ref{NS2})
should be bounded.  Thus,
the tamed system (\ref{NS})-(\ref{NS0}) may serve  as a substitute of (\ref{NS1})-(\ref{NS2}).
Moreover, when we realize Eq. (\ref{NS})-(\ref{NS0}) on a computer,
the value of $N$ can be reset as an arbitrarily large number along with the process of calculations
as long as there is no explosion. So, the term involving $g_N$ serves as some kind of adjustment.

Our main aim in the present paper is to prove the following:
\bt\label{main}
Assume $\nu=1$. Let $\u_0\in\mH^\infty:=\cap_{m\in\mN}\mH^m$
and $[0,\infty)\ni t\mapsto \f(t)\in \mH^m$ be smooth for any $m\in\mN$.
Then there exist a unique smooth velocity field
$$
\u_N\in C^\infty([0,\infty)\times\mR^3;\mR^3)\cap C([0,\infty);\mH^1)
$$
and a pressure function(defined up to a time dependent constant)
$$
p_N\in C^\infty([0,\infty)\times\mR^3;\mR)
$$
solving (\ref{NS})-(\ref{NS0}).
Here the divergence free Sobolev spaces $\mH^m$ are  defined by (\ref{Div}) below.

Moreover, we have the following:
\begin{enumerate}[($1^o$)]
\item
 For any $T>0$
\be
\sup_{t\in[0,T]}\|\u_N(t)\|^2_{\mH^0}
+\int^T_0\|\u_N(s)\|^2_{\mH^1}\dif s\leq
C\left(\|\u_0\|^2_{\mH^0}+\left[\int^T_0\|\f(s)\|_{\mH^0}\dif s\right]^2\right).\label{PI1}
\ee

\item For some absolute constant $C>0$ and any $T,N>0$
\be
\sup_{t\in[0,T]}\|\u_N(t)\|^2_{\mH^1}+\int^T_0\|\u_N(s)\|^2_{\mH^2}\dif s\leq
C_{T,\u_0,\f}\cdot(1+N)\label{PI2},
\ee
\be
\sup_{t\in[0,T]}\|\u_N(t)\|^2_{\mH^2}\leq C'_{T,\u_0,\f}+C_{T,\u_0,\f}\cdot(1+N^2),\label{Ess}
\ee
where
\ce
C_{T,\u_0,\f}&:=&C(T+1)\left(\|\u_0\|^2_{\mH^1}+\|\u_0\|^4_{\mH^1}
+\int^T_0\|\f(s)\|_{\mH^0}^2\dif s\right)\\
C'_{T,\u_0,\f}&:=&C\left(\|\u_0\|^6_{\mH^2}+\sup_{t\in[0,T]}\|\f(t)\|^2_{\mH^0}
+\int^T_0\|\p_s\f(s)\|^2_{\mH^0}\dif s\right).
\de

\item Let $(\u,p)$ be any Leray-Hopf solution of (\ref{NS1})-(\ref{NS2}) on $[0,T]$. If either
$$
\u\in L^\infty([0,T];L^3(\mR^3;\mR^3))
$$
or
$$
\u\in L^q([0,T];L^r(\mR^3;\mR^3)),\quad \frac{3}{r}+\frac{2}{q}=1
$$
holds, then for any bounded domain $\Omega$
\ce
\lim_{N\rightarrow\infty}\int^T_0\int_\Omega|\u_N(t,x)-\u(t,x)|^2\dif x\dif t=0.
\de
\end{enumerate}
\et

\br
By an inequality of Xie \cite{Xie}(see  \cite[Theorem 2.1]{He2}) we know that
\be
\sup_{x\in\mR^3}|\u_N(t,x)|^2\leq \frac{1}{2\pi}\|\Delta\u_N(t)\|_{\mH^0}
\cdot\|\nabla\u_N(t)\|_{\mH^0}.\label{In}
\ee
Hence (\ref{PI2}) and (\ref{Ess}) imply that if $\u_0$ has small enough $\mH^3$-norm and if
$\f$ and $\p_s f$ have small enough $\mH^0$-norm, then
\ce
\sup_{t\in[0,T]}\sup_{x\in\mR^3}|\u_N(t,x)|\leq \sqrt{N},
\de
and (\ref{NS})-(\ref{NS0}) reduces to (\ref{NS1})-(\ref{NS2}), that is,
$\u_N$ solves (\ref{NS1})-(\ref{NS2}).
Moreover, define for any $T, N>0$
\ce
A_{T,N}:=\left\{t\in[0,T]: \sup_{x\in\mR^3}|\u_N(t,x)|\geq \sqrt{N}\right\}.
\de
By (\ref{In}) (\ref{PI1}) and (\ref{PI2}), we have
\ce
\lambda(A_{T,N})&\leq&
\frac{1}{N}\int^T_0\sup_{x\in\mR^3}|\u_N(t,x)|^2\dif t\\
&\leq&\frac{1}{2\pi N}\int^T_0\|\Delta\u_N(t)\|_{\mH^0}
\cdot\|\nabla\u_N(t)\|_{\mH^0}\dif t\\
&\leq&\frac{1}{2\pi N}\left(\int^T_0\|\Delta\u_N(t)\|^2_{\mH^0}\dif t\right)^{1/2}
\left(\int^T_0\|\nabla\u_N(t)\|^2_{\mH^0}\dif t\right)^{1/2}\\
&\leq&\frac{C_{T,\u_0,\f}\cdot(1+N)^{1/2}}
{2\pi N},
\de
where $\lambda(A_{T,N})$ denotes the Lebesgue measure of $A_{T,N}$.
In particular,
$$
\lim_{N\rightarrow\infty}\lambda(A_{T,N})=0,
$$
which shows that as $N$ tends to infinity,
the solution of (\ref{NS})-(\ref{NS0}) satisfies (\ref{NS1})-(\ref{NS2}) at ``almost all''  times.
Note that the following interpolation inequality(see (\ref{Sob}) below)
$$
\sup_{x\in\mR^3}|\u_N(t,x)|^2\leq C\|\u_N(t)\|_{\mH^2}^{3/2}
\cdot\|\u_N(t)\|_{\mH^0}^{1/2}.
$$
For solving (\ref{NS1})-(\ref{NS2}),
an open question is to find an $\a<4/3$ such that for any $N>0$
\ce
\sup_{t\in[0,T]}\|\u_N(t)\|^2_{\mH^2}\leq C_{T,\u_0,f}(1+N^\a).
\de
\er
\br
Define
\ce
\tau_N:=\inf\left\{t>0: \sup_{x\in\mR^3}|\u_N(t,x)|\geq N\right\}.
\de
Clearly, $\u_N(t\wedge\tau_N,x)$ together with some pressure function $p_N$
solves (\ref{NS1})-(\ref{NS2}).
By the uniqueness of solutions to  (\ref{NS})-(\ref{NS0}), we have
$$
\u_{N+1}(t,x)=\u_N(t,x),\quad \forall (t,x)\in [0,\tau_N)\times\mR^3
$$
and
\ce
\tau_N\leq\tau_{N+1}.
\de

Define
\ce
\tau_\infty:=\lim_{N\rightarrow\infty}\tau_N,
\de
and for all $(t,x)\in [0,\tau_N)\times\mR^3$
\ce
\u(t,x):=\u_N(t,x).
\de
Then
$\u$ together with some pressure function $p$ satisfies
(\ref{NS1})-(\ref{NS2}) on $[0,\tau_\infty)\times\mR^3$. In particular, $\tau_\infty$
is the first epoch of irregularity of $\u$(see \cite{Le, Ga}).
\er

\br
Our tamed scheme also works for the periodic case as can be seen from the  proof below.
\er

The approach of proving this theorem is the classical Galerkin's approximation. Instead of working on
$\mH^0$, we shall take $\mH^1$ as our basic space. The special form of the
taming function $g_N$ plays a crucial role. As explained above,
(\ref{NS})-(\ref{NS0}) may serve as an approximation of (\ref{NS1})-(\ref{NS2}).
In fact, in Section 4 we shall  present a new construction for the ``suitable weak solution''
introduced in
\cite{Sc} and \cite{Ca-Ko-Ni} by using $\u_N$. This notion was used to obtain partial regularity of
such type of solutions(cf. \cite{Ca-Ko-Ni} \cite{Li}). We note that it is not at all clear
if weak solutions obtained by the
well known Galerkin approximation procedure(see \cite{Te}) are ``suitable weak solution''.

This paper is organized as follows: in Section 2, some preliminary lemmas are proved.
In Section 3, we prove the existence and uniqueness of weak solution to system (\ref{NS})-(\ref{NS0}).
In Section 4, a suitable weak solution of (\ref{NS1})-(\ref{NS2}) is constructed.
Finally, in Section 5 our main Theorem \ref{main} is proved.

\section{Preliminaries}

Let $C^\infty_0(\mR^3;\mR^3)$
denote the set of smooth functions from $\mR^3$ to $\mR^3$ with
compact supports.
For $p\geq 1$, let $L^p(\mR^3;\mR^3)$ be the vector valued $L^p$-space in which the norm
is denoted by $\|\cdot\|_{L^p}$.
For $m\in\mN_0:=\mN\cup\{0\}$,
let $W^{m,p}$ be the usual Sobolev space on $\mR^3$ with values in $\mR^3$,
i.e., the closure of
$C^\infty_0(\mR^3;\mR^3)$ with respect to the norm:
\ce
\|\u\|_{m,p}=\left(\int_{\mR^3}|(I-\Delta)^{m/2}\u|^p\dif x\right)^{1/p}.
\de
Here as usual $(I-\Delta)^{m/2}$ is defined by Fourier's transformation. This norm is equivalent to the norm
given by
$$
\|\u\|_{m,p}:=\sum_{j=0}^m\left(\int_{\mR^3}|\nabla^j \u|^p\dif x\right)^{1/p},
$$
where $\nabla^j \u$ denotes the total derivative of $\u$ of order $j$.

The following Sobolev type interpolation inequality will be used frequently
and play an essential role in the study of Navier-Stokes equation(cf. \cite[Theorem 2.15]{Ta}).
Let $p,q,r\geq 1$ and $0\leq j<m$. If $m-j-3/p$ is not a non-negative integer and
$$
\frac{1}{r}=\frac{j}{3}+\a\left(\frac{1}{p}-\frac{m}{3}\right)+\frac{1-\a}{q},\quad
\frac{j}{m}\leq\a\leq 1,
$$
then for any $\u\in W^{m,p}\cap L^q(\mR^3;\mR^3)$
\be
\|\nabla^j\u\|_{L^r}\leq C_{m,j,p,q,r}\|\u\|^\a_{m,p}\|\u\|_{L^q}^{1-\a}.\label{Sob}
\ee

Set
\be
\mH^m:=\{\u\in W^{m,2}: \div(\u)=0\},\label{Div}
\ee
where $\div$ is taken in the sense of Schwartz distributions.
Then $(\mH^m,\|\cdot\|_{m,2})$
is a  separable Hilbert space. We shall denote the norm
$\|\cdot\|_{m,2}$ in $\mH^m$ by $\|\cdot\|_{\mH^m}$.
We remark that $\mH^0$ is a closed linear subspace of the Hilbert space $L^2(\mR^3;\mR^3)$
and that by \cite[Theorem 2.21]{Ta} we have both $\mH^1\subset L^r(\mR^3,\mR^3)$ and
$\mH^2\subset W^{1,r}\cap L^q(\mR^3,\mR^3)$ continuously if $r\in[2,6]$, $q\in[2,\infty)$.

Let $R_j$ be the $j$-th Riesz transform(cf. \cite{St}), i.e., for $f\in L^1(\mR^3)$
\ce
R_j(f)(x):=\lim_{\epsilon\rightarrow 0}c_3\int_{|y|\geq\epsilon}\frac{y_j}{|y|^4}f(x-y)\dif y,
\quad j=1,2,3,
\de
where $c_3:=\Gamma(2)/\pi^2$.
Let $P$ be the orthogonal projection operator from $L^2(\mR^3;\mR^3)$ to $\mH^0$.
Then $P$ can be expressed by $R_j$ as(cf.  \cite{Fa-Jo-Ri} \cite{Lions})
\ce
P(\f)=\f+\sum^3_{j=1}\bR R_j(f^j),
\de
where $\bR:=R_1\otimes R_2\otimes R_3$.
It is well known
that $P$ can be restricted to a bounded linear operator from $W^{m,2}$ to $\mH^m$, and that $P$ commutes with the
derivative operators. For any $\u\in\mH^0$ and $\v\in L^2(\mR^3;\mR^3)$, we have
\ce
\<\u,\v\>_{\mH^0}:=\<\u,P\v\>_{\mH^0}=\<\u,\v\>_{L^2}.
\de

Let $\cV$ be defined by
\ce
\cV:=\{\u: \u\in C^\infty_0(\mR^3;\mR^3), \div(\u)=0\}.
\de
We have the following density result.
\bl\label{Le3}
$\cV$ is dense in $\mH^m$ for any $m\in\mN$.
\el
\begin{proof}
For $\rho>0$, let $\om_\rho$ be a mollifier, i.e., positive smooth function on $\mR^3$ with
support in $\{x\in\mR^3: |x|\leq \rho\}$ and $\int_{\mR^3}\om_\rho(x)\dif x=1$.
For $\u\in\mH^m$,
set $\u_\rho(x):=\int_{\mR^3}\u(y)\om_\rho(x-y)\dif y$, then
$\u_\rho\in\cap_{k\in\mN}\mH^k$. For proving the
density of $\cV$ in $\mH^m$, it suffices to prove that $\<\u,\v\>_{\mH^m}=0$ for any $\v\in\cV$ implies
that $\u=0$. Noting that $\v_\rho\in\cV$ if $\v\in\cV$, we have
\ce
0=\<\u,\v_\rho\>_{\mH^m}&=&\<\u,(I-\Delta)^m\v_\rho\>_{\mH^0}\\
&=&\<\u,((I-\Delta)^m\v)_\rho\>_{\mH^0}\\
&=&\<\u_\rho,(I-\Delta)^m\v\>_{\mH^0}\\
&=&\<(I-\Delta)^m\u_\rho,\v\>_{\mH^0}.
\de
Hence, by \cite[Proposition 1.1]{Te} there exists a function $p$ such that
\ce
(I-\Delta)^m\u_\rho=\nabla p,
\de
which leads to
\ce
\int_{\mR^3}|(I-\Delta)^m\u_\rho|^2\dif x=\int_{\mR^3}(I-\Delta)^m\u_\rho\cdot \nabla p~\dif x=0.
\de
Therefore, $\u_\rho=0$ for any $\rho>0$. By taking the limits $\rho\downarrow 0$, we obtain $\u=0$.
The proof is thus complete.
\end{proof}

In the following, for the sake of simplicity, we assume that $\nu=1$.
For  $\u\in\mH^1$ and $\v\in\cV$, define
\ce
\lb A(\u),\v\rb_0:=\<\u,\Delta\v\>_{\mH^0}-\<(\u\cdot\nabla)\u,\v\>_{\mH^0}-\<g_N(|\u|^2)\u,\v\>_{\mH^0},
\de
and
\be
\lb A(\u),\v\rb&:=&\lb A(\u),(I-\Delta)\v\rb_0\label{Es7}\\
&=:&A_1(\u,\v)+A_2(\u,\v)+A_3(\u,\v),\no
\ee
where
\ce
A_1(\u,\v)&:=&\<\u,\Delta \v\>_{\mH^1}\\
A_2(\u,\v)&:=&-\<(\u\cdot\nabla)\u, (I-\Delta)\v\>_{\mH^0}\\
A_3(\u,\v)&:=&-\<g_N(|\u|^2)\u, (I-\Delta)\v\>_{\mH^0}.
\de

Obviously, if $\u\in\mH^2$ the functional $\v\mapsto \lb A(\u),\v\rb_0$ can be identified with the function
$$
\Delta \u-(\u\cdot\nabla)\u-g_N(|\u|^2)\u
$$
which is in $\mH^0$ since by (\ref{Sob}) for $\u\in\mH^2$, $\u\in L^r(\mR^3,\mR^3)$ for all $r\geq 2$
and $\nabla\u\in L^6(\mR^3,\mR^3)$.

We now prepare two lemmas for later use.

\bl\label{Le1}
For any $\u\in\mH^1$ and $\v\in\cV$
\ce
|\lb A(\u),\v\rb|\leq C(1+\|\u\|^3_{\mH^1})\|\v\|_{\mH^3}.
\de
So, by Lemma \ref{Le3}, $\lb A(\u),\cdot\rb$ can be considered as an element in $(\mH^3)'$ with norm bounded
by $C(1+\|\u\|^3_{\mH^1})$. Furthermore, if $\u\in\mH^1$, then
\be
\lb A(\u),\u\rb_0=-\|\nabla\u\|^2_{\mH^0}-\|\sqrt{g_N(|\u|^2)}\cdot|\u|\|^2_{L^2},\label{Es44}
\ee
and if $\u\in\mH^2$, then
\be
\lb A(\u),\u\rb\leq -\frac{1}{2}\| \u\|^2_{\mH^2}-\frac{1}{2}\||\u|\cdot|\nabla\u|\|^2_{L^2}+
(\frac{3}{2}+N)\|\nabla\u\|^2_{\mH^0}+\|\u\|^2_{\mH^0}.\label{Es4}
\ee
\el
\begin{proof}
First of all, it is clear that
\ce
A_1(\u,\v)=\<(I-\Delta)^{1/2}\u, (I-\Delta)^{1/2}\Delta \v\>_{\mH^0}\leq C\|\u\|_{\mH^1}\|\v\|_{\mH^3}.
\de
By the Sobolev  inequality (\ref{Sob}), we have
\ce
A_2(\u,\v)=\<\u^*\cdot\u,\nabla(I-\Delta)\v\>_{\mH^0}
\leq\|\u^*\cdot\u\|_{\mH^0}\cdot\|\nabla(I-\Delta)\v\|_{\mH^0}
\leq C\|\u\|^2_{\mH^1}\cdot\|\v\|_{\mH^3},
\de
where $\u^*$ denotes the transposition of the row vector $\u$,
and
\ce
A_3(\u,\v)\leq \|\u\|^3_{L^6}\cdot\|(I-\Delta) \v\|_{\mH^0}
\leq C\|\u\|^3_{\mH^1}\cdot\|\v\|_{\mH^2}.
\de
Hence
\ce
|\lb A(\u),\v\rb|\leq C(\|\u\|_{\mH^1}+\|\u\|^2_{\mH^1}+\|\u\|^3_{\mH^1})\|\v\|_{\mH^3},
\de
which gives the first assertion.

The equality (\ref{Es44}) clearly follows from
$$
\<(\u\cdot\nabla)\u,\u\>_{\mH^0}=0.
$$

For the inequality (\ref{Es4}), we have
\ce
A_1(\u,\u)&=&-\|(I-\Delta)\u\|^2_{\mH^0}+\<\u, (I-\Delta)\u\>_{\mH^0}\\
&=&-\|\u\|^2_{\mH^2}+\|\nabla\u\|^2_{\mH^0}+\|\u\|^2_{\mH^0},
\de
and by Young's inequality
\ce
A_2(\u,\u)&\leq&\frac{1}{2}\|(I-\Delta)\u\|^2_{\mH^0}
+\frac{1}{2}\|(\u\cdot\nabla)\u\|^2_{\mH^0}\\
&\leq&\frac{1}{2}\|\u\|^2_{\mH^2}+\frac{1}{2}\||\u|\cdot|\nabla\u|\|^2_{\mH^0},
\de
where
\ce
|\u|^2=\sum_{k=1}^3|v^k|^2,\quad |\nabla\u|^2=\sum_{k,i=1}^3|\p_iv^k|^2.
\de

Noting that
\be
g_N(|\u|^2)\geq |\u|^2-(N+1/2)\quad\mbox{ and }\quad g'_N(|\u|^2)\geq 0,\label{Fg}
\ee
we have
\ce
A_3(\u,\u)&=&-\<\nabla (g_N(|\u|^2)\u), \nabla\u\>_{\mH^0}-\<g_N(|\u|^2)\u, \u\>_{\mH^0}\\
&=&-\sum_{k,i=1}^3\int_{\mR^3} \p_i u^k\cdot \p_i (g_N(|\u|^2)u^k)\dif x
-\int_{\mR^3}|\u|^2\cdot g_N(|\u|^2)\dif x\\
&\leq&-\sum_{k,i=1}^3\int_{\mR^3} \p_i u^k\cdot \left(g_N(|\u|^2)\cdot \p_i u^k
-g'_N(|\u|^2)\p_i |\u|^2\cdot u^k\right)\dif x\\\
&=&-\int_{\mR^3}|\nabla \u|^2\cdot g_N(|\u|^2)\dif x
-\frac{1}{2}\int_{\mR^3}g'_N(|\u|^2)|\nabla|\u|^2|^2\dif x\\
&\leq&-\int_{\mR^3}|\nabla \u|^2\cdot |\u|^2\dif x+(N+1/2)\|\nabla\u\|^2_{\mH^0}.
\de
Combining the above calculations yields (\ref{Es4}).
\end{proof}

\bl\label{Le2}
Let $\u_n,\v\in\cV$ and $\u\in\mH^1$.
Let $\Omega:=\mathrm{supp}\v$ and assume that
\ce
\sup_n\|\u_n\|_{\mH^1}<+\infty\mbox{ and }\lim_{n\rightarrow\infty}\|(\u_n-\u)\cdot 1_\Omega\|_{L^2}=0.
\de
Then
\ce
\lim_{n\rightarrow\infty}\lb A(\u_n),\v\rb=\lb A(\u),\v\rb.
\de
\el
\begin{proof}
For $A_1$, we clearly have
\ce
\lim_{n\rightarrow\infty}|A_1(\u_n,\v)-A_1(\u,\v)|
=\lim_{n\rightarrow\infty}|\<(\u_n-\u)\cdot 1_\Omega, (I-\Delta)\Delta \v\>_{\mH^0}|=0.
\de
For $A_2$, obviously
\ce
&&\lim_{n\rightarrow\infty}|A_2(\u_n,\v)-A_2(\u,\v)|\\
&=&\lim_{n\rightarrow\infty}|\<\u^*_n\cdot\u_n-\u^*\cdot\u, \nabla(I-\Delta) \v\>_{\mH^0}|=0.
\de
For $A_3$, we have because $|g_N'|\leq C$ and by (\ref{Sob})
\ce
&&\lim_{n\rightarrow\infty}|A_3(\u_n,\v)-A_3(\u,\v)|\\
&\leq&C\sup_{x\in\mR^3}|(I-\Delta)\v(x)|\cdot \lim_{n\rightarrow\infty}
(\|\u_n\|^2_{\mH^1}\cdot\|(\u_n-\u)\cdot 1_\Omega\|_{L^2})= 0.
\de
The proof is complete.
\end{proof}

\section{Existence and Uniqueness of Weak Solution}

We first give the following definition of a generalized solution to (\ref{NS})-(\ref{NS0})
according to \cite{La}.
\bd\label{Def}
Let $T>0$,  $\u_0\in\mH^0$ and $\f\in L^2([0,T];\mH^0)$.
A measurable vector field  $\u$ on $[0,T]\times\mR^3$
is called a generalized solution of (\ref{NS})-(\ref{NS0}) if
\begin{enumerate}[($1^o$)]
\item $\u\in L^\infty([0,T]; L^4(\mR^3;\mR^3))\cap L^2([0,T];\mH^1)$;
\item for any $\v\in\cV$ and $t\in[0,T]$,
\be
\<\u(t),{\v}\>_{\mH^0}&=&\<\u_0,{\v}\>_{\mH^0}-\int^t_0\<\nabla\u(s),\nabla \v\>_{\mH^0}\dif s
-\int^t_0\<(\u(s)\cdot\nabla)\u(s),\v\>_{\mH^0}\dif s\no\\
&&-\int^t_0\<g_N(|\u(s)|^2)\u(s),\v\>_{\mH^0}\dif s+\int^t_0\<\f(s),\v\>_{\mH^0}\dif s;\label{Eq}
\ee
\item $\lim_{t\downarrow 0}\|\u(t)-\u_0\|_{L^2}=0$.
\end{enumerate}
\ed
\br
In the definition of a weak solution of the classical Navier-Stokes equation(i.e., $N=\infty$),
$\u\in L^\infty([0,T]; L^4(\mR^3;\mR^3))$ is replaced by $\u\in L^\infty([0,T]; L^2(\mR^3;\mR^3))$.
For the former, we have uniqueness, but no existence, and for the later, we have existence,
but no uniqueness(cf. \cite{Te}).
\er

The following proposition is a well known consequence of Lemma \ref{Le3}(cf. \cite{Ga}).
\bp\label{P1}
 Let $\u$ be  a weak solution of (\ref{NS})-(\ref{NS0}) in the sense of
Definition \ref{Def}. Then we have for any $\v\in C^1([0,T];\mH^1)$ with $\v(T)=0$
\be
&&\int^T_0\<\u(t),\p_t\v(t)\>_{\mH^0}\dif t+\<\u_0,\v(0)\>_{\mH^0}\no\\
&=&\int^T_0\<\nabla\u(t),\nabla \v(t)\>_{\mH^0}\dif t
+\int^T_0\<(\u(t)\cdot\nabla)\u(t),\v(t)\>_{\mH^0}\dif t\no\\
&&+\int^T_0\<g_N(|\u(t)|^2)\u(t),\v(t)\>_{\mH^0}\dif t-\int^T_0\<\f(t),\v(t)\>_{\mH^0}\dif t.\label{Es9}
\ee
Moreover, the following energy equality holds
\be
&&\|\u(t)\|^2_{\mH^0}+2\int^t_0\|\nabla\u(s)\|^2_{\mH^0}\dif s
+2\int^t_0\|\sqrt{g_N(|\u(s)|^2)}|\u(s)|\|^2_{\mH^0}\dif s\no\\
&=&\|\u_0\|^2_{\mH^0}+2\int^t_0\<\f(s),\u(s)\>_{\mH^0}\dif s,\quad\forall t\in[0,T].\label{En}
\ee
\ep
\begin{proof}
We have by Sobolev's inequality (\ref{Sob}) and H\"older's inequality
\ce
&&\int^T_0|\<(\u(t)\cdot\nabla)\u(t),\v(t)\>_{\mH^0}|\dif t\\
&\leq& \int^T_0\|\u(t)\|_{L^4}\|\nabla\u(t)\|_{L^2}\cdot\|\v(t)\|_{L^4}\dif t\\
&\leq&C\cdot \sup_{t\in[0,T]}\|\u(t)\|_{L^4}\cdot \left(\int^T_0\|\u(t)\|^2_{\mH^1}\dif t\right)^{1/2}
\left(\int^T_0\|\v(t)\|^2_{\mH^1}\dif t\right)^{1/2}
\de
and
\ce
&&\int^T_0|\<g_N(|\u(t)|^2)\u(t),\v(t)\>_{\mH^0}|\dif t\\
&\leq&\int^T_0\int_{\mR^3}|\u(t,x)|^3\cdot|\v(t,x)|\dif x\dif t\\
&\leq&\int^T_0\|\u(t)\|^3_{L^4}\cdot \|\v(t)\|_{L^4}\dif t\\
&\leq&C\cdot \sup_{t\in[0,T]}\|\u(t)\|^3_{L^4}\cdot
\int^T_0\|\v(t)\|_{\mH^1}\dif t.
\de
Hence, the right hand side of (\ref{Es9}) is well defined. By using (\ref{Eq}),
equality (\ref{Es9}) and the energy equality (\ref{En}) follow from an appropriate approximation.
\end{proof}
We may now prove the following uniqueness result.
\bt
(Uniqueness) Let $\u$ and $\tilde \u$ be two generalized solutions of (\ref{NS})-(\ref{NS0}) in the sense of
Definition \ref{Def}. Then $\u=\tilde \u$.
\et
\begin{proof}
Set
$$
\w:=\u-\tilde\u.
$$
Then analogously to proving (\ref{En}) we obtain
\ce
\|\w(t)\|^2_{\mH^0}&=&-2\int^t_0\|\nabla\w(s)\|^2_{\mH^0}\dif s
-2\int^t_0\<\w(s),(\u(s)\cdot\nabla)\u(s)-(\tilde\u(s)\cdot\nabla)\tilde\u(s)\>_{\mH^0}\dif s\\
&&-2\int^t_0\<\w(s),g_N(|\u(s)|^2)\u(s)-g_N(|\tilde\u(s)|^2)\tilde\u(s)\>_{\mH^0}\dif s\\
&=:&I_1+I_2+I_3.
\de
By (\ref{Sob}) and Young's inequality, we have for any $\epsilon>0$
\be
&&\|\u^*(s)\cdot\u(s)-\tilde\u^*(s)\cdot\tilde\u(s)\|^2_{\mH^0}\no\\
&\leq&\||\w(s)|(|\u(s)|+|\tilde\u(s)|)\|^2_{\mH^0}\no\\
&\leq&2\|\w(s)\|^2_{L^4}(\|\u(s)\|^2_{L^4}+\|\tilde\u(s)\|^2_{L^4})\no\\
&\leq&2C_{1,0,2,2,4}\|\w(s)\|^{3/2}_{\mH^1}\|\w(s)\|^{1/2}_{\mH^0}(\|\u(s)\|^2_{L^4}+\|\tilde\u(s)\|^2_{L^4})\no\\
&\leq&C_\epsilon M_{\u,\tilde\u}\|\w(s)\|^2_{\mH^0}+\epsilon \|\w(s)\|^2_{\mH^1},\label{Es8}
\ee
where
\ce
M_{\u,\tilde\u}:=\mathrm{ess}\sup_{s\in[0,T]}\left(\|\u(s)\|^8_{L^4}
+\|\tilde \u(s)\|^8_{L^4}\right)<+\infty.
\de

Hence, for $\epsilon=1/4$
\ce
I_2&=&2\int^t_0\<\nabla\w(s),(\u^*(s)\cdot\u(s)-\tilde\u^*(s)\cdot\tilde\u(s))\>_{\mH^0}\dif s\\
&\leq&\frac{1}{2}\int^t_0\|\nabla\w(s)\|^2_{\mH^0}\dif s
+2\int^t_0\|\u^*(s)\cdot\u(s)-\tilde\u^*(s)\cdot\tilde\u(s))\|^2_{\mH^0}\dif s\\
&\leq&\int^t_0\|\w(s)\|^2_{\mH^1}\dif s+2C  M_{\u,\tilde\u}\int^t_0\|\w(s)\|^2_{\mH^0}\dif s.
\de

On the other hand, we also have for any $\epsilon>0$
\ce
&&|\<\w(s),g_N(|\u(s)|^2)\u(s)-g_N(|\tilde\u(s)|^2)\tilde\u(s)\>_{\mH^0}|\\
&\leq&\int_{\mR^3}|\w(s,x)|^2 g_N(|\u(s,x)|^2)\dif x\\
&&+\int_{\mR^3}|\w(s,x)|\cdot|g_N(|\u(s,x)|^2)-g_N(|\tilde\u(s,x)|^2)|\cdot|\tilde\u(s,x)|\dif x\\
&\leq&4\||\w(s)|(|\u(s)|+|\tilde\u(s)|)\|^2_{\mH^0}\\
&\leq&C_\epsilon M_{\u,\tilde\u}\|\w(s)\|^2_{\mH^0}+\epsilon \|\w(s)\|^2_{\mH^1},
\de
where we have used that $|g'_N|\leq C$ and (\ref{Es8}).

Hence for $\epsilon=1$ we obtain by integrating with respect to $\dif s$
\ce
I_3\leq\int^t_0\|\w(s)\|^2_{\mH^1}\dif s+C  M_{\u,\tilde\u}\int^t_0\|\w(s)\|^2_{\mH^0}\dif s.
\de
For $I_1$, we have
\ce
I_1=-2\int^t_0\|\w(s)\|^2_{\mH^1}\dif s+2\int^t_0\|\w(s)\|^2_{\mH^0}\dif s.
\de
Combining all estimates on $I_1,I_2$ and $I_3$, we obtain
\ce
\|\w(t)\|^2_{\mH^0}\leq C  M_{\u,\tilde\u}\int^t_0\|\w(s)\|^2_{\mH^0}\dif s.
\de
So, by Gronwall's inequality
\ce
\w(t)=0,\quad\forall t>0,
\de
and the assertion is proved.
\end{proof}

Let us now prove the following existence result.
\bt\label{Th1}
Let $\u_0\in\mH^1$ and
$\f\in L^2_{\mathbf{loc}}([0,\infty);\mH^0)$. Then there exists at least one solution $\u$
to (\ref{NS})-(\ref{NS0}) in the sense of
Definition \ref{Def}. Moreover,

\begin{enumerate}[($1^\circ$)]
\item $\u\in C([0,\infty);\mH^1)
\cap L^2_{\mathbf{loc}}([0,\infty);\mH^2)$, $\p_t\u\in L^2_{\mathbf{loc}}([0,\infty);\mH^0)$, and
\be
&&\|\u(t)\|^2_{\mH^1}+\int^t_0\left(\|\u(s)\|^2_{\mH^2}
+\||\u(s)|\cdot|\nabla\u(s)|\|^2_{L^2}\right)\dif s\no\\
&\leq& C(1+N)(t+1)\left(\|\u_0\|^2_{\mH^1}
+\int^t_0\|\f(s)\|_{\mH^0}^2\dif s\right),\quad t\geq 0, \label{Es11}
\ee
where $N$ is as in (\ref{NS})-(\ref{NS0}).
\item
\be
\|\u(t)\|_{\mH^0}\leq \|\u_0\|_{\mH^0}
+\int^t_0\|\f(s)\|_{\mH^0}\dif s,\quad t\geq 0\label{OP3}
\ee
and
\be
&&\int^t_0\left(\|\nabla\u(s)\|^2_{\mH^0}
+\||\u(s)|\cdot \sqrt{g_N(|\u(s)|^2)}\|^2_{\mH^0}\right)\dif s\no\\
&\leq&C\left(\|\u_0\|^2_{\mH^0}+\left[\int^t_0\|\f(s)\|_{\mH^0}\dif s\right]^2\right),\quad t\geq 0.\label{Es01}
\ee

\item There exists a real function $p(t,x)$ with $\nabla p\in L^2_{\mathbf{loc}}([0,\infty);L^2(\mR^3;\mR^3))$
such that for almost all $t\geq 0$
$$
\p_t\u(t)=\Delta\u(t)-(\u(t)\cdot\nabla)\u(t)+\nabla p(t)-g_N(|\u(t)|^2)\u(t)+\f(t).
$$
\end{enumerate}

\et
\begin{proof}
We use the standard Galerkin approximation, and only consider the finite time interval $[0,T]$.
Let $\{\e_k,k\in\mN\}\subset\cV$ be a complete orthonormal basis of $\mH^1$ such that
$\mathrm{span}\{e_k,k\in\mN\}$ is dense in $\mH^3$.
Fix $n\in\mN$.
For $\y=(y^1,\cdots,y^n)\in\mR^n$, set
\ce
\y\cdot \e&:=&\sum_{i=1}^n y^i\e_i\in\cV\\
\b_n(\y)&:=&(\lb A(\y\cdot\e),\e_1\rb,\cdots,\lb A(\y\cdot\e),\e_n\rb)\\
\f_n(t)&:=&(\<\rho_n*\f(t),\e_1\>_{\mH^1},\cdots,\<\rho_n*\f(t),\e_n\>_{\mH^1}),
\de
where $\rho_n$ is a family of mollifiers and $\rho_n*f(t)$ denotes the convolution on $\mR^3$.

Consider the following ordinary differential equation
\ce
\dif \y_n(t)/\dif t=\b_n(\y_n(t))+\f_n(t)
\de
subject to the initial condition
\ce
\y_n(0)=(\<\u_0,\e_1\>_{\mH^1},\cdots,\<\u_0,\e_n\>_{\mH^1}).
\de
By Lemma \ref{Le1}, we have for some $C_{n,N}>0$
\ce
\<\y,\b_n(\y)\>_{\mR^n}\leq C_{n,N}|\y|^2.
\de
Moreover, it is easy to see that $\y\mapsto \b_n(\y)$ is a smooth function.
Hence, by the theory of ODE there is a unique $\y_n(t)$ satisfying
\ce
\y_n(t)=\y_n(0)+\int^t_0\b_n(\y_n(s))\dif s+\int^t_0\f_n(s)\dif s,\ \ t\geq 0.
\de

Set
\ce
\u_n(t)&:=&\y_n(t)\cdot \e=\sum_{i=1}^ny^i_n(t)\e_i\\
\Pi_n A(\u_n(t))&:=&\sum_{i=1}^n\lb A(\u_n(t)),\e_i\rb\e_i\\
\Pi_n\f(t)&:=&\sum_{i=1}^n\<\rho_n*\f(t),\e_i\>_{\mH^1}\e_i.
\de
Then
\be
\p_t\u_n(t)=\Pi_n A(\u_n(t))+\Pi_n\f(t),\label{OP4}
\ee
and for any $n\geq k$
\be
\<\u_n(t),\e_k\>_{\mH^1}&=&\<\u_n(0),\e_k\>_{\mH^1}+\int^t_0\<\Pi_n A(\u_n(s)),\e_k\>_{\mH^1}\dif s
+\int^t_0\<\Pi_n\f(s),\e_k\>_{\mH^1}\dif s\no\\
&=&\<\u_0,\e_k\>_{\mH^1}+\int^t_0\lb A(\u_n(s)),\e_k\rb\dif s
+\int^t_0\<\rho_n*\f(s),\e_k\>_{\mH^1}\dif s.\label{Es6}
\ee
We then have by (\ref{Es6}) and (\ref{Es4})
\ce
\|\u_n(t)\|^2_{\mH^1}&=&\|\u_n(0)\|^2_{\mH^1}+2\int^t_0\lb A(\u_n(s)),\u_n(s)\rb\dif s
+2\int^t_0\<\rho_n*\f(s),\u_n(s)\>_{\mH^1}\dif s\\
&\leq&\|\u_0\|^2_{\mH^1}-\int^t_0\|\u_n(s)\|^2_{\mH^2}\dif s
-\int^t_0\||\u_n(s)|\cdot|\nabla\u_n(s)|\|^2_{L^2}\dif s\\
&&+(3+2N)\int^t_0\|\nabla\u_n(s)\|^2_{\mH^0}\dif s+2\int^t_0\|\u_n(s)\|^2_{\mH^0}\dif s\\
&&+2\int^t_0\|\rho_n*\f(s)\|_{\mH^0}\|\u_n(s)\|_{\mH^2}\dif s\\
&\leq&\|\u_0\|^2_{\mH^1}-\frac{1}{2}\int^t_0\|\u_n(s)\|^2_{\mH^2}\dif s
-\int^t_0\||\u_n(s)|\cdot|\nabla\u_n(s)|\|^2_{L^2}\dif s\\
&&+C_N\int^t_0\|\u_n(s)\|^2_{\mH^1}\dif s+2\int^t_0\|\f(s)\|^2_{\mH^0}\dif s,
\de
which implies by Gronwall's inequality that for all $n\in\mN$
\be
\sup_{t\in[0,T]}\|\u_n(t)\|^2_{\mH^1}+\int^T_0\|\u_n(s)\|^2_{\mH^2}\dif s
\leq C_{N,T,\f}.\label{Es1}
\ee

For fixed $k\in\mN$, set
$$
G^{(k)}_n(t):=\<\u_n(t),\e_k\>_{\mH^1}.
$$
(\ref{Es1}) implies that $G^{(k)}_n, n\in\mN$, are uniformly bounded. Furthermore,
for any $0\leq r<t\leq T$, we have from (\ref{Es6}) and Lemma \ref{Le1}
\ce
|G^{(k)}_n(t)-G^{(k)}_n(r)|&\leq&\int^t_r|\lb A(\u_n(s)),\e_k\rb|\dif s
+\int^t_r|\<\rho_n*\f(s),\e_k\>_{\mH^1}|\dif s\\
&\leq&C\|\e_k\|_{\mH^3}\int^t_r(1+\|\u_n(s)\|^3_{\mH^1})\dif s
+\|\e_k\|_{\mH^2}\int^t_r\|\f(s)\|_{\mH^0}\dif s,
\de
that $G^{(k)}_n, n\in\mN$, are also equi-continuous by (\ref{Es1}).

Thus, by Ascoli-Arzel\`a lemma, we may substract a subsequence $n_l$
such that $G^{(k)}_{n_l}(t)$ uniformly converges to a continuous function $G^{(k)}(t)$.
This subsequence may depend on $k$. However, by a diagonalization method, we may find a common
subsequence (still denoted by $n$) such that for any $k\in\mN$
$$
\lim_{n\rightarrow\infty}\sup_{t\in[0,T]}|G^{(k)}_n(t)-G^{(k)}(t)|=0.
$$
Thus, by (\ref{Es1}) and since closed balls of $\mH^1$ are weakly compact, there exits a
$\u\in L^\infty([0,T];\mH^1)$ such that $t\mapsto\u(t)$ is weakly continuous, that is,
continuous from $[0,T]$ to $\mH^1$ equipped with the weak topology, and
for any $\v\in\mH^1$
\ce
\lim_{n\rightarrow\infty}\sup_{t\in[0,T]}|\<\u_n(t)-\u(t),\v\>_{\mH^1}|=0.
\de
In particular, for any $\tilde\v\in\mH^0$, taking $\v=(I-\Delta)^{-1}\tilde\v\in\mH^2$ we get
\ce
\lim_{n\rightarrow\infty}\sup_{t\in[0,T]}|\<\u_n(t)-\u(t),\tilde\v\>_{\mH^0}|=0.
\de
By Helmholtz-Weyl orthogonal decomposition(cf. \cite{Te}),
we further have for any $\w\in L^2(\mR^3;\mR^3)$  and $T>0$
\be
\lim_{n\rightarrow\infty}\sup_{t\in[0,T]}|\<\u_n(t)-\u(t),\w\>_{L^2}|=0.\label{Es33}
\ee

Moreover, by (\ref{Es1}) and  (\ref{Es33}), we also have
\ce
\int^T_0\|\u(s)\|^2_{\mH^2}\dif s\leq \liminf_{n\rightarrow\infty}\int^T_0\|\u_n(s)\|^2_{\mH^2}\dif s
<+\infty.
\de

Let us now check that $\u$ constructed above is a solution of (\ref{NS})-(\ref{NS0}).

Let $\cO\subset\mR^3$ be any bounded  domain. Using (\ref{Es1}) (\ref{Es33}) and
the following Friedrichs' inequality(see \cite[p.176]{La}):  For any
$\epsilon>0$, there exist  $N_\epsilon\in\mN$
functions $h_i^\epsilon\in L^2(\cO)$ such that for any $w\in W^{1,2}_0(\cO)$
\be
\int_\cO|w(x)|^2\dif x\leq
\sum^{N_\epsilon}_{i=1}\left(\int_\cO w(x) h_i^\epsilon(x)\dif x\right)^2
+\epsilon\int_\cO|\nabla w(x)|^2\dif x,\label{Fre}
\ee
we obviously have that for any $\cO'\subset\overline\cO'\subset\cO$
\be
\lim_{n\rightarrow\infty}\sup_{t\in[0,T]}\int_{\cO'}|\u_n(t,x)-\u(t,x)|^2\dif x=0.\label{Es222}
\ee

For any $k\in\mN$, let the support of $\e_k$ be contained in some bounded  domain $\cO_k$.
By Lemma \ref{Le2}, (\ref{Es222}) and the dominated convergence theorem, we have
\ce
\int^t_0\lb  A(\u_n(s)),\e_k\rb\dif s\rightarrow \int^t_0\lb  A(\u(s)),\e_k\rb\dif s,
~ \mbox{ as}~ n\rightarrow\infty.
\de
Now taking the limits $n\rightarrow\infty$ for (\ref{Es6}) we obtain that for any $k\in\mN$
and $t\geq 0$
\ce
\<\u(t),\e_k\>_{\mH^1}=\<\u(0),\e_k\>_{\mH^1}+\int^t_0\lb A(\u(s)),\e_k\rb\dif s
+\int^t_0\<\f(s),(I-\Delta)\e_k\>_{\mH^0}\dif s.
\de
By an easy approximation and Lemma \ref{Le1}, we further have for any $\v\in\mH^3$ and $t\geq 0$
\ce
\<\u(t),\v\>_{\mH^1}=\<\u(0),\v\>_{\mH^1}+\int^t_0\lb A(\u(s)),\v\rb\dif s
+\int^t_0\<\f(s),(I-\Delta)\v\>_{\mH^0}\dif s.
\de
Since $(I-\Delta)^{-1}\v\in\mH^3$ for any $\v\in\mH^3$, from the above equality and (\ref{Es7})
we obtain (\ref{Eq}).

Let us now prove ($1^o$) and ($2^o$).
Since $\u\in L^\infty([0,T];\mH^1)\cap L^2([0,T];\mH^2)$ for any $T>0$,
by (\ref{Eq}) we have for almost all $t\geq 0$
\be
\p_t\u(t)=\Delta\u(t)-P((\u(t)\cdot\nabla)\u(t))-P(g_N(|\u(t)|^2)\u(t))+\f(t) \mbox{ in $\mH^0$}\label{Op5}
\ee
and $\p_t\u(t)\in L^2([0,T];\mH^0)$.
($3^o$) now follows from  (\ref{Eq}) and \cite[Proposition 1.1]{Te}.

Taking inner products with $\u(t)$ in $\mH^0$ for both sides of (\ref{Op5}),
we have by (\ref{Es44})
\ce
\frac{\dif \|\u(t)\|^2_{\mH^0}}{\dif t}&=&2\<A(\u(t)),\u(t)\>_{\mH^0}+2\<\f(t),\u(t)\>_{\mH^0}\no\\
&\leq&-2\|\nabla\u(t)\|^2_{\mH^0}-2\||\u(t)|\cdot \sqrt{g_N(|\u(t)|^2)}\|^2_{\mH^0}\no\\
&&+2\|\f(t)\|_{\mH^0}\cdot\|\u(t)\|_{\mH^0},
\de
and thus
\ce
\frac{\dif \|\u(t)\|_{\mH^0}}{\dif t}\leq\|\f(t)\|_{\mH^0}.
\de
Both together implies ($2^o$).

Since $\mH^2\subset\mH^1\subset\mH^0$ forms an evolutional triple,
we also have by \cite[p. 260, Lemma 1.2]{Te} and (\ref{Op5})
\be
\|\u(t)\|^2_{\mH^1}=\|\u(0)\|^2_{\mH^1}+2\int^t_0\lb A(\u(s)),\u(s)\rb\dif s
+2\int^t_0\<\f(s),\u(s)\>_{\mH^1}\dif s.\label{PP1}
\ee
Since the right hand side of (\ref{PP1}) is continuous in $t$ and we have already know that
$t\mapsto \u(t)\in\mH^1$ is weakly continuous, it follows that this map is even (strongly) continuous in $\mH^1$.
Furthermore, the right hand side of (\ref{PP1}) is by  (\ref{Es4}) and (\ref{Es01}) bounded by
\ce
&&\|\u_0\|^2_{\mH^1}-\int^t_0\|\u(s)\|^2_{\mH^2}\dif s
-\int^t_0\||\u(s)|\cdot|\nabla\u(s)|\|^2_{L^2}\dif s\\
&&+(3+2N)\int^t_0\|\nabla\u(s)\|^2_{\mH^0}\dif s+2\int^t_0\|\u(s)\|^2_{\mH^0}\dif s\\
&&+2\int^t_0\|\f(s)\|_{\mH^0}\|\u(s)\|_{\mH^2}\dif s\\
&\leq&-\frac{1}{2}\int^t_0\|\u(s)\|^2_{\mH^2}\dif s
-\int^t_0\||\u(s)|\cdot|\nabla\u(s)|\|^2_{L^2}\dif s\\
&&+C(1+N)(t+1)\left(\|\u_0\|^2_{\mH^1}+\int^t_0\|\f(s)\|^2_{\mH^0}\dif s\right),
\de
which gives (\ref{Es11}).
The proof is thus finished.
\end{proof}

\section{Existence of Suitable Weak Solution}

In this section we shall use the result of previous section to
give a new construction the ``suitable weak solution'' to (\ref{NS1})-(\ref{NS2})
that was introduced in \cite{Sc}  and \cite{Ca-Ko-Ni}.
\bt\label{main2}
Given $T>0$, let $\u_0\in\mH^0$ and
$\f\in L^2([0,T];\mH^0)$. There exists a Leray-Hopf weak solution $(\u,p)$ to (\ref{NS1})-(\ref{NS2})
satisfying the following generalized energy inequality:

If $\phi\in C^\infty_0((0,T)\times\mR^3)$ is nonnegative, then
\be
&&2\int^T_0\int_{\mR^3}|\nabla\u(t,x)|^2\phi(t,x)\dif x\dif s\no\\
&\leq&\int^T_0\int_{\mR^3}\Big[|\u(t,x)|^2(\p_t\phi(t,x)+\Delta\phi(t,x))+2(\u(t,x)\cdot\f(t,x))\phi(t,x)\no\\
&&+(|\u(t,x)|^2+2p(t,x))\<\u(t,x),\nabla\phi(t,x)\>_{\mR^3}\Big]\dif x\dif t.\label{Gen}
\ee
\et
\begin{proof}
Choose $\u^N_0\in\mH^1$ being such that
\ce
\lim_{N\rightarrow\infty}\|\u^N_0-\u_0\|_{\mH^0}=0.
\de
Let $\u_N$ be the unique solution of (\ref{NS})-(\ref{NS0}) constructed in Theorem \ref{Th1}
with initial value $\u^N_0$. By ($2^o$) of Theorem \ref{Th1}, we have
\be
\sup_{t\in[0,T]}\|\u_N(t)\|^2_{\mH^0}+\int^T_0\left(\|\u_N(s)\|^2_{\mH^1}
+\|\sqrt{g_N(|\u_N(s)|^2)}\cdot|\u_N(s)|\|^2_{L^2}\right)\dif s\leq C_{T,\u_0,\f}.\label{Op3}
\ee
Let $q\in[2,\infty)$ and $r\in(2,6]$
\ce
\frac{3}{r}+\frac{2}{q}=\frac{3}{2}.
\de
By (\ref{Sob}) and (\ref{Op3}) we have
\be
\int^T_0\|\u_N(t)\|^q_{L^r}\dif t\leq C^q_{1,0,2,2,r}
\int^T_0\|\u_N(t)\|^2_{\mH^1}\cdot\|\u_N(t)\|^{q-2}_{\mH^0}\dif t\leq C_{T,\u_0,\f,r,q}.\label{Op11}
\ee

As in the proof of Theorem \ref{Th1}, one may find a subsequence still denoted by $N$
and $\u\in L^\infty([0,T];\mH^0)\cap L^2([0,T];\mH^1)$ such that
for any $\w\in L^2(\mR^3;\mR^3)$
\be
\lim_{N\rightarrow\infty}\sup_{t\in[0,T]}|\<\u_N(t)-\u(t),\w\>_{L^2}|=0.\label{Es30}
\ee
We shall now prove that in fact for any bounded domain $\Omega\subset\mR^3$
\be
\lim_{N\rightarrow\infty}\int^T_0\int_\Omega|\u_N(t,x)-\u(t,x)|^2\dif x\dif t=0.\label{Lim}
\ee
Now, let $\Omega\subset\overline\Omega\subset\Omega'$ and $\rho\geq 0$ be a smooth cutoff function
with $\rho(x)=1$ on $\Omega$ and $\rho(x)=0$ on $\mR^3-\Omega'$. Then by Friedrichs' inequality (\ref{Fre})
we have
\ce
&&\int^T_0\int_\Omega|\u_N(t,x)-\u(t,x)|^2\dif x\dif t\\
&\leq&\int^T_0\int_{\Omega'}|\u_N(t,x)-\u(t,x)|^2\cdot\rho^2(x)\dif x\dif t\\
&\leq&\sum^{N_\epsilon}_{i=1}\int^T_0\left(\int_{\Omega'} (\u_N(t,x)-\u(t,x))\rho(x) h_i(x)\dif x\right)^2\dif t\\
&&+\epsilon\int^T_0\int_{\Omega'}|\nabla((\u_N(t)-\u(t))\rho)(x)|^2\dif x\dif t\\
&=:&I_1(k,\epsilon)+I_2(k,\epsilon).
\de
The second term is bounded by
\ce
I_2(k,\epsilon)\leq\epsilon\cdot C_\rho\int^T_0
\left(\|\u_N(t)\|^2_{\mH^1}+\|\u(t)\|^2_{\mH^1}\right)\dif t\leq C_{\rho,T,\u_0,\f}
\cdot\epsilon.
\de
By (\ref{Es30}) we have for the first term
\ce
\lim_{N\rightarrow \infty}I_1(k,\epsilon)\leq T\sum^{N_\epsilon}_{i=1}
\lim_{N\rightarrow \infty}\sup_{t\in[0,T]}\left|\int_{\mR^3}
 (\u_N(t,x)-\u(t,x))\rho(x) h_i^\epsilon(x) 1_{\Omega'}(x)\dif x\right|^2=0.
\de
Hence (\ref{Lim}) follows by the arbitrariness of $\epsilon$.

As in the proof of Theorem \ref{Th1}, in order to verify that $\u$ is a Leray-Hopf weak solution,
it suffices to prove that for any $\v\in\cV$
\ce
\lim_{N\rightarrow\infty}\int^t_0\<g_{N}(|\u_N(s)|^2)\u_N(s), \v\>_{\mH^0}\dif s=0.
\de
This is true, because by (\ref{Op11}), this limit is dominated by
\ce
&&\|\v\|_{L^\infty}\cdot\limsup_{N\rightarrow\infty}\int^t_0\int_{\mR^3}|\u_N(s,x)|^3\cdot
1_{\{|\u_N(s,x)|^2\geq N\}}\dif x\dif s\\
&\leq&\|\v\|_{L^\infty}\cdot\limsup_{N\rightarrow\infty}\left(\int^t_0\|\u_N(s)\|^{10/3}_{L^{10/3}}\dif s\right)^{9/10}
\cdot\left(\int^t_0\int_{\mR^3}1_{\{|\u_N(s,x)|^2\geq N\}}\dif x\dif s\right)^{1/10}\\
&\leq&C_{\v,T,\u_0,\f}\cdot\limsup_{N\rightarrow\infty}\left(\frac{1}{N}\int^t_0
\|\u_N(s)\|^2_{\mH^0}\dif s\right)^{1/10}=0.
\de

Let us now check the generalized energy inequality.
Note that by \cite[Proposition 1.1]{Te},
for some $p_N\in L^2([0,T]; L^2_{\mathbf{loc}}(\mR^3,\mR^3))$ with $\nabla p_N\in L^2([0,T];L^2(\mR^3;\mR^3))$
\be
\p_t\u_N=\Delta\u_N(t)-(\u_N(t)\cdot\nabla)\u_N(t)+\nabla p_N(t)-g_N(|\u_N(t)|^2)\u_N(t)+\f(t).\label{Op2}
\ee
Taking inner product with $2\u_N\phi$ and integrating by parts yield
\be
&&2\int^T_0\int_{\mR^3}|\nabla\u_N|^2\phi~\dif x\dif s
+2\int^T_0\int_{\mR^3}g_N(|\u_N|^2)|\u_N|^2\phi~\dif x\dif s\no\\
&=&\int^T_0\int_{\mR^3}\left[|\u_N|^2(\p_t\phi+\Delta\phi)
+2(\u_N\cdot\f)\phi+(|\u_N|^2-2p_N)\<\u_N,\nabla\phi\>_{\mR^3}\right]\dif x\dif t.\label{Gen1}
\ee
Since $\u\mapsto\int^T_0\int_{\mR^3}|\nabla\u|^2\phi~\dif x\dif s$ is lower-semi continuous in
$L^2([0,T],\mH^0)$, the limit of the
left hand side of (\ref{Gen1}) is greater than the left hand side of (\ref{Gen}).
Since $\phi$ has compact support, by (\ref{Op11})(with $q=r=10/3$) and (\ref{Lim}),
each term of the right hand side of (\ref{Gen1}) converges to
the corresponding term of the right hand side of (\ref{Gen})
except for the term including $p_N$. For this term, one takes the divergence
(in the sense of Schwartz distributions) for (\ref{Op2})
and then gets
\be
\Delta p_N=\div((\u_N(t)\cdot\nabla)\u_N(t)+g_N(|\u_N(t)|^2)\u_N(t)).\label{be1}
\ee
Noting that
\ce
g^{9/8}_N(r)\cdot r^{9/16}\leq g_N(r)\cdot  r,
\de
we have by (\ref{Op3})
\be
&&\int^T_0\int_{\mR^3}|g_N(|\u_N(t,x)|^2)\cdot \u_N(t,x)|^{9/8}\dif x\dif t\no\\
&\leq&\int^T_0\int_{\mR^3}g_N(|\u_N(t,x)|^2)\cdot |\u_N(t,x)|^2\dif x\dif t\leq C_{T,\u_0,\f}.\label{be2}
\ee

Moreover, by H\"older's inequality and (\ref{Op11})
\be
&&\int^T_0\int_{\mR^3}|(\u_N(t)\cdot\nabla)\u_N(t)|^{9/8}\dif x\dif t\no\\
&\leq&\left(\int^T_0\|\u_N(t)\|^{18/7}_{L^{18/7}}\dif t\right)^{7/16}
\cdot\left(\int^T_0\|\u_N(t)\|^2_{\mH^1}\dif t\right)^{9/16}\leq C_{T,\u_0,\f}.\label{be3}
\ee
By (\ref{Sob}), we know that
\be
\int^T_0\|p_N(t)\|^{9/8}_{L^{9/5}}\dif t\leq
C^{9/8}_{1,0,9/8,9/8,9/5}\int^T_0\|p_N(t)\|^{9/8}_{1,9/8}\dif t.\label{be4}
\ee
Since the operator $\nabla(-\Delta)^{1/2}$ is bounded in $L^p(\mR^3)$ for all $p>1$(cf. \cite{St}),
it follows from (\ref{be1})-(\ref{be3}) that the right hand side of (\ref{be4}) is uniformly bounded in $N$.

Therefore, there exists a subsequence(still denoted by $N$) and
$$
p\in L^{9/8}([0,T];L^{9/5}(\mR^3;\mR^3))
$$
such that
\be
p_N\rightarrow p\mbox{ weakly in }L^{9/8}([0,T];L^{9/5}(\mR^3;\mR^3)).\label{PP1a}
\ee
On the other hand, by  (\ref{Op11}) with $q=12, r=9/4$ again we have
\ce
\int^T_0\|\u_N(t)\|^{12}_{L^{9/4}}\dif t\leq C_{T,\u_0,\f}.
\de
So, by (\ref{Lim}) and (\ref{PP1a})
\ce
&&\lim_{N\rightarrow\infty}\left|\int^T_0\int_{\mR^3}\<p_N\u_N-p\u,\nabla\phi\>_{\mR^3}(t,x)\dif x\dif t\right|\\
&\leq&\lim_{N\rightarrow\infty}\left|\int^T_0\int_{\mR^3}(p_N-p)(t,x)\cdot\<\u,\nabla\phi\>_{\mR^3}(t,x)\dif x\dif t\right|\\
&&+\lim_{N\rightarrow\infty}\left|\int^T_0\int_{\mR^3}p_N(t,x)\cdot\<\u_N-\u,\nabla\phi\>_{\mR^3}(t,x)\dif x\dif t\right|\\
&\leq&\lim_{N\rightarrow\infty}\left(\int^T_0\|p_N(t)\|^{9/8}_{9/5}\dif t\right)^{8/9}
\left(\int^T_0\||\u_N(t)-\u(t)|\cdot|\nabla\phi(t)|\|^9_{9/4}\dif t\right)^{1/9}=0.
\de
The proof is complete.
\end{proof}

\section{Proof of Theorem \ref{main}}

For $T>0$, we write $S_T:=[0,T]\times\mR^3$, and
\ce
\|\u\|_{L^q(S_T)}:=\left\{
\begin{array}{ll}
\left(\int_{S_T}|\u(s,x)|^q\dif x\dif s\right)^{1/q},\ \ \ q\in[1,\infty), \\
\sup_{(s,x)\in S_T}|\u(s,x)|, \ \ \ q=\infty.
\end{array}
\right.
\de
Let $T_t$ be the Gaussian heat semigroup on $\mR^3$, i.e., for $h\in L^1(\mR^3)$
\ce
T_th(x):=\frac{1}{(4\pi t)^{3/2}}\int_{\mR^3}e^{-|x-y|^2/(4t)}h(y)\dif y.
\de
Define
$$
B(\u,\u)(t,x):=\sum^3_{i=1}\int^t_0\p_i T_{t-s}\Big[u^i(s)\u(s)-\sum^3_{j=1}\bR R_ju^j(s)u^i(s)\Big](x)\dif s.
$$
Then, the unique solution $\u$ constructed in Theorem \ref{Th1} for fixed $N$
satisfies the following integral equation(cf. \cite[Theorem 4.4]{Fa-Jo-Ri}):
\ce
\u(t,x)=\f_N(t,x)-B(\u,\u)(t,x),
\de
where  $\f_N$ is defined by
\ce
\f_N(t,x):=T_t\u_0(x)-\int^t_0T_{t-s}P(g_N(|\u(s)|^2)\u(s)-\f(s))(x)\dif s.
\de

Let us recall the following regularity result due to Fabes-Jones-Riviere \cite[Theorem 3.4]{Fa-Jo-Ri}.
\bt\label{regu}
Assume $\u\in L^q(S_T)$ for some $q\geq 5$. Let $k\in\mN$ be such that $k+1<q$. If
$$
D^\alpha_x D^j_t\f_N\in L^{q/(|\a|+2j+1)}(S_T) \quad \mbox{ whenever $|\a|+2j\leq k$},
$$
then also
$$
D^\alpha_x D^j_t\u\in L^{q/(|\a|+2j+1)}(S_T) \quad \mbox{ whenever $|\a|+2j\leq k$}.
$$
Here $D^\a$ denotes the usual derivative operator and $\a=(\a_1,\a_2,\a_3)$ denotes a multi index,
$|\a|=\a_1+\a_2+\a_3$.
\et

We are now in a position to prove our main Theorem \ref{main}.

\vspace{3mm}

{\it Proof of Theorem \ref{main}}:

\vspace{3mm}

Let us now use  induction  to prove that for $k\in\mN$
\ce
\u\in L^{10\cdot(5/3)^{k-1}}(S_T),\quad D^\alpha_x D^j_t\u\in L^{2\cdot(5/3)^k}(S_T) \ \
\mbox{ for $|\a|+2j\leq 2$},\quad (\sP_k).
\de
First of all, by the Sobolev inequality (\ref{Sob}) and (\ref{Es11}), we have
\ce
\|\u\|_{L^{10}(S_T)}^{10}&\leq& C^{10}_{2,0,2,6,10}\int^T_0\|\u(s)\|^2_{\mH^2}\cdot\|\u(s)\|^8_{L^6}\dif s\\
&\leq& C\sup_{s\in[0,T]}\|\u(s)\|^8_{\mH^1}\int^T_0\|\u(s)\|^2_{\mH^2}\dif s<+\infty.
\de
Thus
$$
g_N(|\u|^2)\u\in L^{10/3}(S_T).
$$
Since $\u_0$ and $\f$ are smooth, by the $L^p$-theory of singular integrals of parabolic type
(cf. \cite{La-So-Ur}),
we know
$$
D^\alpha_x D^j_t\f_N\in L^{10/3}(S_T) \quad \mbox{ for $|\a|+2j\leq 2$}.
$$
By Theorem \ref{regu}, we get
$$
D^\alpha_x D^j_t\u\in L^{10/3}(S_T) \quad \mbox{ for $|\a|+2j\leq 2$}.
$$
Hence $(\sP_1)$ holds.

Now suppose that $(\sP_k)$ holds.
Using the Sobolev inequality (\ref{Sob})  and (\ref{Es11}) again, we have
\ce
\|\u\|_{L^{10\cdot(5/3)^{k}}(S_T)}^{10\cdot(5/3)^{k}}
&\leq& C^{10\cdot(5/3)^{k}}_{2,0,2\cdot(5/3)^k,6,10\cdot(5/3)^{k}}
\int^T_0\|\u(s)\|^{2\cdot(5/3)^k}_{2,2\cdot(5/3)^k}\cdot\|\u(s)\|^{8\cdot(5/3)^{k}}_{L^6}\dif s\\
&\leq& C\sup_{s\in[0,T]}\|\u(s)\|^{8\cdot(5/3)^{k}}_{\mH^1}
\int^T_0\|\u(s)\|^{2\cdot(5/3)^k}_{2,2\cdot(5/3)^k}\dif s<+\infty.
\de
Thus,
$$
g_N(|\u|^2)\u\in L^{2\cdot(5/3)^{k+1}}(S_T)
$$
and
$$
D^\alpha_x D^j_t\f_N\in L^{2\cdot(5/3)^{k+1}}(S_T) \quad \mbox{ for $|\a|+2j\leq 2$}.
$$
By Theorem \ref{regu} again, we get
$$
D^\alpha_x D^j_t\u\in L^{2\cdot(5/3)^{k+1}}(S_T) \quad \mbox{ for $|\a|+2j\leq 2$},
$$
i.e., $(\sP_{k+1})$ holds.

Therefore,
$$
D^\alpha_x D^j_t\u\in \cap_{q>1}L^{q}(S_T) \quad \mbox{ for $|\a|+2j\leq 2$},
$$

Using Theorem \ref{regu}, by the induction method, one finds that for any $m\in\mN$
$$
D^\alpha_x D^j_t\u\in \cap_{q>1}L^{q}(S_T) \quad \mbox{ for $|\a|+2j\leq m$}.
$$
The smoothness of $\u$ now follows from the Sobolev embedding theorem
(cf. \cite{La-So-Ur}). In particular, for any $k,m\in\mN$ and $T>0$
$$
\p^k_t\u\in L^2([0,T];\mH^m)
$$
and for all $t\geq 0$ and $m\in\mN$
\be
\p_t\u(t)=\Delta\u(t)-P((\u(t)\cdot\nabla)\u(t))-P(g_N(|\u(t)|^2)\u(t))+\f(t) \mbox{ in $\mH^m$}.\label{Op9}
\ee

We now prove the estimate (\ref{Ess}).
First of all, taking inner products with $\p_t\u(t)$ in $\mH^0$ for both sides of (\ref{Op9}),
we have by (\ref{Es7}), (\ref{OP4}) and Young's inequality
\ce
\|\p_t\u(t)\|^2_{\mH^0}&=&\<\Delta\u(t), \p_t\u(t)\>_{\mH^0}
-\<(\u(t)\cdot\nabla)\u(t),\p_t\u(t)\>_{\mH^0}\\
&&-\<g_N(|\u(t)|^2)\u(t),\p_t\u(t)\>_{\mH^0}+\<\f(t),\p_t\u(t)\>_{\mH^0}\\
&\leq&-\frac{1}{2}\cdot\frac{\dif\|\nabla\u(t)\|^2_{\mH^0}}{\dif t}
+\|\p_t\u(t)\|_{\mH^0}\|(\u(t)\cdot\nabla)\u(t)\|_{L^2}\\
&&-\frac{1}{2}\cdot\frac{\dif\|\tilde g_N(|\u(t)|^2)\|_{L^1}}{\dif t}
+\frac{1}{4}\|\p_t\u(t)\|^2_{\mH^0}+\|\f(t)\|^2_{\mH^0}\\
&\leq&-\frac{1}{2}\cdot\frac{\dif\|\nabla\u(t)\|^2_{\mH^0}}{\dif t}
+\frac{1}{2}\|\p_t\u(t)\|^2_{\mH^0}+\||\u(t)|\cdot|\nabla\u(t)|\|^2_{L^2}\\
&&-\frac{1}{2}\cdot\frac{\dif\|\tilde g_N(|\u(t)|^2)\|_{L^1}}{\dif t}+\|\f(t)\|^2_{\mH^0},
\de
where
$$
0\leq \tilde g_N(r):=\int^r_0 g_N(s)\dif s\leq \frac{r^2}{2}.
$$

Integrating both sides from $0$ to $T$ with respect to $t$
and using (\ref{Es11}) yields that
\be
\frac{1}{2}\int^T_0\|\p_t\u(t)\|^2_{\mH^0}\dif t&\leq&-\frac{1}{2}\left(\|\nabla\u(T)\|^2_{\mH^0}
-\|\nabla\u(0)\|^2_{\mH^0}\right)\no\\
&&+\int^T_0\left(\||\u(t)|\cdot|\nabla\u(t)|\|^2_{L^2}
+\|\f(t)\|^2_{\mH^0}\right)\dif t\no\\
&&-\frac{1}{2}\left(\|\tilde g_N(|\u(T)|^2)\|_{L^1}-\|\tilde g_N(|\u(0)|^2)\|_{L^1}\right)\no\\
&\leq&\frac{1}{2}\|\u(0)\|^2_{\mH^1}+\frac{1}{4}\|\u(0)\|^4_{L^4}\no\\
&&+C(1+N)(1+T)\left(\|\u_0\|^2_{\mH^1}+\int^T_0\|\f(s)\|_{\mH^0}^2\dif s\right)\no\\
&\leq& C(1+N)(1+T)\left(\|\u_0\|^2_{\mH^1}+\|\u_0\|^4_{\mH^1}
+\int^T_0\|\f(s)\|_{\mH^0}^2\dif s\right).\label{Op1}
\ee

On the other hand,
differentiating (\ref{Op9}) with respect to $t$ and then
taking the inner products with $\p_t\u(t)$ in $\mH^0$ gives that
\ce
\frac{1}{2}\frac{\dif \|\p_t\u(t)\|^2_{\mH^0}}{\dif t}
&=&\<\Delta\p_t\u(t), \p_t\u(t)\>_{\mH^0}
-\<(\p_t\u(t)\cdot\nabla)\u(t),\p_t\u(t)\>_{\mH^0}\\
&&-\<\p_t(g_N(|\u(t)|^2)\u(t)),\p_t\u(t)\>_{\mH^0}+\<\p_t\f(t),\p_t\u(t)\>_{\mH^0}\\
&=&-\|\nabla\p_t\u(t)\|^2_{\mH^0}
-\<\nabla\p_t\u(t), \u^*(t)\cdot\p_t\u(t)\>_{\mH^0}\\
&&-\<g_N(|\u(t)|^2)\p_t\u(t),\p_t\u(t)\>_{\mH^0}+\<\p_t\f(t),\p_t\u(t)\>_{\mH^0}\\
&&-\<g'_N(|\u(t)|^2)\p_t|\u(t)|^2\u(t),\p_t\u(t)\>_{\mH^0}\\
&\leq&-\frac{1}{2}\|\nabla\p_t\u(t)\|^2_{\mH^0}
+\frac{1}{2}\||\u(t)|\cdot|\p_t\u(t)|\|^2_{\mH^0}\\
&&-\|\sqrt{g_N(|\u(t)|^2)}\cdot|\p_t\u(t)|\|^2_{\mH^0}
+\frac{1}{2}\|\p_t\f(t)\|^2_{\mH^0}+\frac{1}{2}\|\p_t\u(t)\|^2_{\mH^0}\\
&&-\frac{1}{2}\|\sqrt{g'_N(|\u(t)|^2)}\cdot|\p_t|\u(t)|^2|\|^2_{L^2}\\
&\leq&-\frac{1}{2}\|\nabla\p_t\u(t)\|^2_{\mH^0}
-\frac{1}{2}\||\u(t)|\cdot|\p_t\u(t)|\|^2_{\mH^0}\\
&&+(N+1)\|\p_t\u(t)\|^2_{\mH^0}+\frac{1}{2}\|\p_t\f(t)\|^2_{\mH^0},
\de
where we have used that $g_N(|\u|^2)\geq |\u|^2-(N+1/2)$.

Integrating both sides and using (\ref{Op1}) yield that for any $t\geq 0$
\be
\|\p_t\u(t)\|^2_{\mH^0}&\leq&\|\p_t\u(0)\|^2_{\mH^0}+2(N+1)\int^t_0\|\p_s\u(s)\|^2_{\mH^0}\dif s+
\int^t_0\|\p_s\f(s)\|^2_{\mH^0}\dif s\no\\
&\leq& C(1+N^2)(t+1)\left(\|\u_0\|^2_{\mH^1}+\|\u_0\|^4_{\mH^1}
+\int^t_0\|\f(s)\|_{\mH^0}^2\dif s\right)\no\\
&&+C(\|\u(0)\|^6_{\mH^2}+\|\f(0)\|^2_{\mH^0})+\int^t_0\|\p_s\f(s)\|^2_{\mH^0}\dif s,\label{Es90}
\ee
where we have used that $\|\p_t\u(0)\|^2_{\mH^0}\leq C(\|\u(0)\|^6_{\mH^2}+\|\f(0)\|^2_{\mH^0})$ by (\ref{Op9}).

Lastly,  taking inner products with $\u(t)$ in $\mH^1$ for both sides of (\ref{Op9}), by (\ref{Es4}) we get
\ce
\<\p_t\u(t),\u(t)\>_{\mH^1}&=&\lb A(\u(t)),\u(t)\rb+\<\f(t),\u(t)\>_{\mH^1}\\
&\leq&-\frac{1}{4}\|\u(t)\|^2_{\mH^2}-\frac{1}{2}\||\u(t)|\cdot|\nabla\u(t)|\|^2_{L^2}\\
&&+(\frac{3}{2}+N)\|\nabla\u(t)\|^2_{\mH^0}+\|\u(t)\|^2_{\mH^0}+\|\f(t)\|^2_{\mH^0}.
\de
This implies that
\ce
\sup_{t\in[0,T]}\|\u(t)\|^2_{\mH^2}&\leq& 4(\frac{3}{2}+N)\sup_{t\in[0,T]}\|\nabla\u(t)\|^2_{\mH^0}
+4\sup_{t\in[0,T]}\|\u(t)\|^2_{\mH^0}\\
&&+4\sup_{t\in[0,T]}\|\f(t)\|^2_{\mH^0}+4\sup_{t\in[0,T]}|\<\p_t\u(t),\u(t)\>_{\mH^1}|\\
&\leq& C(1+N^2)(T+1)\left(\|\u_0\|^2_{\mH^1}+\|\u_0\|^4_{\mH^1}
+\int^T_0\|\f(s)\|_{\mH^0}^2\dif s\right)\\
&&+4\sup_{t\in[0,T]}\|\f(t)\|^2_{\mH^0}
+8\sup_{t\in[0,T]}\|\p_t\u(t)\|^{2}_{\mH^0}+\frac{1}{2}\sup_{t\in[0,T]}\|\u(t)\|^2_{\mH^2},
\de
where we used (\ref{Es11}) in the second step.

So, by (\ref{Es90})
\ce
\sup_{t\in[0,T]}\|\u(t)\|^2_{\mH^2}&\leq& C(1+N^2)(T+1)\left(\|\u_0\|^2_{\mH^1}+\|\u_0\|^4_{\mH^1}
+\int^T_0\|\f(s)\|_{\mH^0}^2\dif s\right)\\
&&+C\left(\|\u_0\|^6_{\mH^2}+\sup_{t\in[0,T]}\|\f(t)\|^2_{\mH^0}+\int^T_0\|\p_s\f(s)\|^2_{\mH^0}\dif s\right),
\de
which leads to the desired estimate (\ref{Ess}).
Since (\ref{PI1}) follows from (\ref{Op3}) (\ref{Es01}),
and (\ref{PI2}) follows from (\ref{Es11}), ($2^o$) is proved.

For ($3^o$), suppose that the conclusion is false. Then, there exist a bounded domain
$\Omega$, $\epsilon>0$ and a subsequence $N_k$ such that
\ce
\int^T_0\int_\Omega|\u_{N_k}(t,x)-\u(t,x)|^2\dif x\dif t\geq \epsilon.
\de

Starting from this subsequence and using Theorem \ref{main2}, one may find a suitable solution $\tilde\u$
and another subsequence $N'_k$ of $N_k$ such that(see (\ref{Lim}))
\ce
\lim_{k\rightarrow\infty}\int^T_0\int_\Omega|\u_{N'_k}(t,x)-\tilde\u(t,x)|^2\dif x\dif t=0.
\de
Under the assumptions of ($3^o$), by  well known results
(cf. \cite[Theorem 4.2 and Theorem 7.2]{Ga}) we know that
\ce
\tilde\u=\u,~ a.e.
\de
This contradiction proves ($3^o$). \ \ \ \ \ \ \ \ \ \ \ \ \ \ \ \ \ \
\ \ \ \ \ \ \ \ \ \ \ \ \ \ \ \ \ \ \ \ \ \ \ \ \ \ \ \ \ \ \ \ \ \ \ \ \ \ \ \ \ \ \ \ \ \ \ \ \ \ \ \ \ \  $\square$

\end{document}